\documentclass[12pt,a4paper]{article}

\usepackage[left=2cm,right=2cm,    top=2cm,bottom=2cm,bindingoffset=0cm]{geometry}

\usepackage{subfig}
\usepackage{flushend}
\usepackage{indentfirst}
\usepackage{graphics}
\usepackage{amsmath}
\usepackage{epstopdf}
\usepackage{float}
\usepackage{amssymb}
\usepackage[font=footnotesize,labelfont=bf]{caption}
\usepackage[labelsep=period]{caption}

\newcommand{\me}{\mathbf}
\newcommand{\mr}{\mathbb}
\newcommand{\mt}{\mathsf}
\newcommand{\md}{\mathcal}
\newcommand{\eps}{\varepsilon}
\newcommand{\ld}{\left}
\newcommand{\rd}{\right}

\newcommand{\be}{\begin{equation}}
\newcommand{\ee}{\end{equation}}

\newtheorem{thm}{Theorem}[section]
\newtheorem{lem}{Lemma}[section]

\newtheorem{ozn}{Definition}[section]

\begin{document}

\title{\textbf{Filtering Problem for Random Processes with Stationary Increments }}

\date{}

\maketitle

\noindent Columbia International Publishing\\
Contemporary Mathematics and Statistics
(2015) Vol. 3 No. 1 pp. 8-27\\
doi:10.7726/cms.2015.1002

\vspace{20pt}

\author{\textbf{Maksym Luz}, \textbf{Mykhailo Moklyachuk}$^*$, \\\\
 {Department of Probability Theory, Statistics and Actuarial
Mathematics, \\
Taras Shevchenko National University of Kyiv, Kyiv 01601, Ukraine}\\
$^{*}$Corresponding Author: Moklyachuk@gmail.com}\\\\\\

\noindent \textbf{Abstract.} \hspace{2pt}  This paper deals with the  problem  of optimal mean-square filtering of the linear functionals $A{\xi}=\int_{0}^{\infty}a(t)\xi(-t)dt$ and $A_T{\xi}=\int_{0}^Ta(t)\xi(-t)dt$ which depend on the
unknown values of random process $\xi(t)$ with stationary $n$th
increments from observations of  process $\xi(t)+\eta(t)$ at points
$t\leq0$, where $\eta(t)$ is a stationary process uncorrelated with $\xi(t)$. We propose the values of mean-square errors and spectral characteristics of optimal linear
estimates of the functionals when spectral
densities of the processes are known. In the case where we can operate only with a set of admissible spectral densities
relations that determine the
least favorable spectral densities and the minimax spectral
characteristics are proposed.\\

\noindent \textbf{Keywords:} \hspace{2pt} Random process with $n$th stationary
increments, minimax-robust estimate, mean-square error, least favorable spectral density, minimax spectral
characteristic.

\section{Introduction}

The problems of estimating unobserved values of random processes play an important role in both theoretical and applied probability. Widely developed methods of extrapolation, interpolation and filtering  of stationary random processes and sequences were proposed by Kolmogorov
 \cite{Kolmogorov} and Wiener
\cite{Wiener}. Since that time a lot of generalizations of stationary processes have been established; some of which are discussed, for example, in  books by Yaglom \cite{Yaglom:1987a,Yaglom:1987b}. One important generalization was proposed by Yaglom in  paper \cite{Yaglom:1955}, in which he considered random processes with stationary $n$th increments. He described spectral representation of stationary increments and canonical factorization of spectral density, solved extrapolation problem for unknown value of random process with stationary increments and discussed a few examples for given spectral densities. The further investigations of such random processes were made in  articles by Pinsker \cite{Pinsk:1955}, Yaglom and Pinsker \cite{Pinsk:1954}, Yaglom \cite{Yaglom:1957}.

The crucial assumption of most of the papers dedicated to one of problems of estimating values of random processes is that spectral densities of involved random processes are exactly known. However, the established results can't be directly applied to  practical issues, because we operate with estimates of spectral densities, not with exact views. Often instead of estimates researchers operate with sets of admissible spectral densities. In this case they can use minimax-robust approach to the estimation problem. This method provides us an estimate that minimizes
the maximum of mean-square errors for all spectral densities from the
 set of admissible spectral
densities simultaneously.  Grenander \cite{Grenander} was the first one who proposed this method for solving extrapolation problem for stationary processes. Later
Franke
\cite{Franke:1985}, Franke and Poor \cite{Franke} applied minimax-robust method to extrapolation and interpolation problems. They used convex optimization for investigations. In the papers by Moklyachuk \cite{Moklyachuk:1989}
-\cite{Moklyachuk:2008} minimax-robust extrapolation, interpolation and filtering  problems are studied.   The papers by Moklyachuk and Masyutka
\cite{Moklyachuk:2006} -\cite{Moklyachuk:2012} are dedicated to minimax-robust extrapolation, interpolation and filtering  problems for vector-valued stationary processes and sequences.
 Dubovetska, Moklyachuk and Masyutka \cite {Dubovetska1} solved a problem of minimax-robust interpolation of another generalization of stationary processes --  periodically correlated sequences. In the papers  \cite {Dubovetska4} -\cite{Dubovetska8} Dubovetska and Moklyachuk solved minimax-robust extrapolation, interpolation and filtering  problems for periodically correlated processes.

 Minimax-robust extrapolation, interpolation and filtering  problems for stochastic  sequences and random processes with $n$th stationary increments  are solved by Luz and Moklyachuk in \cite{Luz1} -- \cite{Luz6}. In particular,  minimax-robust filtering problem for such sequences is considerer in papers \cite{Luz3,Luz5}. Random processes with stationary increments with continuous time are considered in the article \cite{Luz4} by Luz and Moklyachuk, where authors investigate extrapolation problem.

In this article we propose solution to the problem of optimal mean-square estimation of the linear functionals $A_T{\xi}=\int_{0}^Ta(t)\xi(-t)dt$ and $A{\xi}=\int_{0}^{\infty}a(t)\xi(-t)dt$ of unknown values of a random process $\xi(t)$ with stationary $n$th increments  from observation of a process $\xi(t)+\eta(t)$ at points $t\leq0$, where a process $\eta(t)$ is stationary uncorrelated with $\xi(t)$. we propose formulas for calculating the value of mean-square errors and the spectral characteristics of the optimal estimates of the functionals $A_T{\xi}$ and $A{\xi}$  under the condition that spectral densities of the processes are exactly known. In the case when spectral densities are not known, but  sets of admissible spectral densities are available, relations which determine  least favorable densities and  minimax-robust spectral characteristics for different classes of spectral densities are found.

\section{\Large{Stationary random increment process. Spectral representation}}

In this section we present a brief review of the spectral theory of random processes with stationary increment developed in \cite{Yaglom:1955}.

\begin{ozn}
For a given random process $\{\xi(t),t\in \mr R\}$ a process
\begin{equation}
\label{oznachPryrostu_cont}
\xi^{(n)}(t,\tau)=(1-B_{\tau})^n\xi(t)=\sum_{l=0}^n(-1)^lC_n^l\xi(t-l\tau),
\end{equation}
where $B_{\tau}$ is a backward shift operator with step $\tau\in  \mr R$, such that
$B_{\tau}\xi(t)=\xi(t-\tau)$, is called the random $n$th increment with step  $\tau\in \mr R$.
\end{ozn}

Random $n$th increment process $\xi^{(n)}(t,\tau)$ satisfies  relations:
\begin{equation}
\xi^{(n)}(t,-\tau)=(-1)^n\xi^{(n)}(t+n\tau,\tau), \label{tot1_cont}
\end{equation}
\begin{equation} and
\label{tot2_cont}
\xi^{(n)}(t,k\tau)=\sum\nolimits_{l=0}^{(k-1)n}A_l\xi^{(n)}(t-l\tau,\tau),\quad
\forall k\in  \mr N,
\end{equation}
where coefficients $\{A_l,l=0,1,2,\ldots,(k-1)n\}$ are taken from a representation
\[(1+x+\ldots+x^{k-1})^n=\sum_{l=0}^{(k-1)n
}A_lx^l.
\]

\begin{ozn}
Random $n$th increment process $\xi^{(n)}(t,\tau)$ generated
by a random process $\{\xi(t),t\in  R\}$ is wide sense
stationary if the mathematical expectations
\[
\mathsf E\xi^{(n)}(t_0,\tau)=c^{(n)}(\tau),
\]
\[
{\mathsf E}
\xi^{(n)}(t_0+t,\tau_1)\xi^{(n)}(t_0,\tau_2)=D^{(n)}(t,\tau_1,\tau_2)
\]
exist for all $t_0,\tau,t,\tau_1,\tau_2$ and do not depend on $t_0$. The function $c^{(n)}(\tau)$ is called the mean value of the $n$th increment and the function
$D^{(n)}(t,\tau_1,\tau_2)$ is called the structural function of the stationary $n$th increment (or the structural function of $n$th
order of the random process $\{\xi(t),t\in \mr R\}$).

The random process $\{\xi(t),t\in \mr R\}$ which determines the
stationary $n$th increment process $\xi^{(n)}(t,\tau)$ by formula
(\ref{oznachPryrostu_cont})  is called the process with stationary $n$th
increments.
\end{ozn}

The following theorem provodes representations of mean value and structural function of stationary random $n$th increment.

\begin{thm}  \label{thm1_cont}
The mean value $c^{(n)}(\tau)$ and the structural function
$D^{(n)}(m,\tau_1,\tau_2)$ of a random stationary $n$th increment
process  $\xi^{(n)}(t,\tau)$ can be represented in the following
forms
\begin{equation}
\label{serFnaR_cont}
c^{(n)}(\tau)=c\tau^n,
\end{equation}
 \begin{equation}\label{strFnaR_cont}D^{(n)}(t,\tau_1,\tau_2)=\int_{-\infty}^{\infty} e^{i\lambda
  t } (1-e^{-i\tau_1\lambda})^n(1-e^{i\tau_2\lambda})^n\frac{(1+ \lambda^2)^{n}}
{\lambda^{2n}}dF(\lambda),
\end{equation}
where  $c$ is a constant, $F(\lambda)$ is a left-continuous
nondecreasing bounded function with $F(-\infty)=0$. The constant $c$
and the function $F(\lambda)$ are determined uniquely by the
increment process $\xi^{(n)}(t,\tau)$.

On the other hand, the function $c^{(n)}(\tau)$ which has the form
$(\ref{serFnaR_cont})$ with a constant $c$ and the function
$D^{(n)}(t,\tau_1,\tau_2)$ which has the form $(\ref{strFnaR_cont})$ with a
function $F(\lambda)$ which satisfies the indicated conditions are
the mean value and the structural function of some stationary $n$th
increment process $\xi^{(n)}(t,\tau)$.
\end{thm}

From representation $(\ref{strFnaR_cont})$ of the structural function
 of the stationary
$n$th increment process $\xi^{(n)}(t,\tau)$ and the Karhunen theorem
(see Karhunen $\cite{Karhunen}$), we obtaine the following spectral representation of
the stationary $n$th increment process
$\xi^{(n)}(t,\tau)$:
\begin{equation}
\label{predZnaR_cont} \xi^{(n)}(t,\tau)=\int_{-\infty}^{\infty}
e^{it \lambda }(1-e^{-i\lambda\tau})^n\frac{(1+i\lambda )^n}{(i\lambda)^n}dZ(\lambda),
\end{equation}
where $Z(\lambda)$ -- is a random process with independent increments on $\mr R$ connected with \emph{a spectral function} $F(\lambda)$ from representation ($\ref{strFnaR_cont}$) by the relation
\begin{equation}
\label{FtaZ_cont}
 \mathsf E|Z(t_2)-Z(t_1)|^2=F(t_2)-F(t_1)<\infty\quad \mbox{for all } t_2>t_1.
\end{equation}

Without loss of generality we will consider random increments $\xi^{(n)}(t,\tau)$ with step  $\tau>0$ and mean value $0$.

\section{Filtering problem}
Consider a random process $\{\xi(t),t\in\mr R\}$ which defines a stationary $n$th increment process $\xi^{(n)}(t,\tau)$ with  an absolutely continuous spectral function  $F(\lambda)$ which has a spectral density $f(\lambda)$. Assume that we observe another process $\zeta(t)=\xi(t)+\eta(t)$ on the time interval $t\leq0$. Filtering problem means that we need to restore values of the original process $\xi(t)$ at points $t\leq0$.

In this article we focus on finding mean-square optimal estimates of linear functionals $A_T\xi=\int_0^T a(t)\xi(-t)dt$ and $A\xi=\int_0^{\infty} a(t)\xi(-t)dt$ under the assumption that a noise process $\eta(t)$ is an uncorrelated with $\xi(t)$ stationary process with spectral density $g(\lambda)$.

 For the further discussions we  require the spectral densities $f(\lambda)$ and $g(\lambda)$ of the random processes $\xi(t)$ and $\eta(t)$ satisfy minimality conditions
\be \label{int_cont_umova1}
 \int_{-\infty}^{\infty}\frac{|\gamma(\lambda)|^2}{f(\lambda)+\frac{\lambda^{2n}} {(1+\lambda^2 )^n}g(\lambda)}d\lambda<\infty,\quad\int_{-\infty}^{\infty}
\frac{|\gamma(\lambda)|^2\lambda^{2n}}{|1-e^{i\lambda
\tau}|^{2n}(1+\lambda^2)^{n}(f(\lambda)+\frac{\lambda^{2n}} {(1+\lambda^2 )^n}g(\lambda))} d\lambda<\infty\ee for some non-zero function of the exponential type
$\gamma(\lambda)=\int_{0}^{\infty}\alpha(t)e^{i\lambda t}dt$. We  also require a function $\me a_{\tau}(t)$, $t\geq-\tau n$, being defined later satisfy conditions
\be\label{cond_fro_function_a_f_cont_st_n}\int_0^{\infty}|\me a_{\tau}(t-\tau n)|dt<\infty,\quad \int_0^{\infty}t|\me a_{\tau}(t-\tau n)|^2dt.\ee

 To solve the filtering problem for the functional $A\xi$ and the random processes $\xi(t)$ and $\eta(t)$ described above we represent the functional $A\xi=\int_0^{\infty} a(t)\xi(-t)dt$ in the form
\[A\xi=A\zeta-A\eta,\] where $A\zeta=\int_0^{\infty}a(t)\zeta(-t)dt$, $A\eta=\int_0^{\infty}a(t)\eta(-t)dt$. As we can find the exact value of the functional $A\zeta$ from the observations $\zeta(t)$ at points $t\leq0$, it's sufficient to find an estimate $\widehat{A}\eta$ of the functional $A\eta$ in order to build the estimate $\widehat{A}\xi:$
\be\label{mainformula}\widehat{A}\xi=A\zeta-\widehat{A}\eta.\ee
Moreover, the values of mean-square errors $\Delta(f,g;\widehat{A}\xi)=\mt E |A\xi-\widehat{A}\xi|^2$ and $\Delta(f,g;\widehat{A}\eta)=\mt E|A\eta-\widehat{A}\eta|^2$ of the estimates $\widehat{A}\xi$ and $\widehat{A}\eta$ respectively satisfy the equalities
\[\Delta\ld(f,g;\widehat{A}\xi\rd)=\mt E \ld|A\xi-\widehat{A}\xi\rd|^2= \mt
E\ld|A\zeta-A\eta-A\zeta+\widehat{A}\eta\rd|^2=\mt
E\ld|A\eta-\widehat{A}\eta\rd|^2=\Delta\ld(f,g,\widehat{A}\eta\rd).\]

In this paper we do not try to find the optimal estimate $\widehat{A}\eta$ by direct minimizing  mean-square error as a function of admissible estimates. We just use  properties of minimal distances in Hilbert space $H=L_2(\Omega, \mathfrak{F},
\mt P)$ of random variables of the second order. In this space the available observations $\{\xi^{(n)}(t,\tau)+\eta^{(n)}(t,\tau):t\leq0\}$, $\tau>0$, generate  the closed
linear subspace $H^{0}(\xi^{(n)}_{\tau}+\eta^{(n)}_{\tau})$ of the Hilbert space $H=L_2(\Omega, \mathfrak{F},
\mt P)$. Thus the optimal estimate $\widehat{A}\eta$ is a projection of the element $A\eta$ of the space $H$ on the subspace $H^{0}(\xi^{(n)}_{\tau}+\eta^{(n)}_{\tau})$ and the value of mean-square error $\Delta(f,g;\widehat{A}\eta)$ is the minimal distance from the element $A\eta$ from the space $H$ to the subspace $H^{0}(\xi^{(n)}_{\tau}+\eta^{(n)}_{\tau})$. These properties give us two conditions characterizing the optimal estimate $\widehat{A}\eta$:

1) $ \widehat{A}\eta\in H^{0}\ld(\xi^{(n)}_{\tau}+\eta^{(n)}_{\tau}\rd) $;

2) $\ld({A}\eta-\widehat{A}\eta\rd)\perp
H^{0}\ld(\xi^{(n)}_{\tau}+\eta^{(n)}_{\tau}\rd)$.

The obtained conditions allow us to find a spectral characteristic $h_{\tau}(\lambda)$ of the estimate $\widehat{A}\eta$. To be able to use them we need to describe the spectral representations of the involved processes.

The
stationary random process $\eta(t)$ with spectral function $G(\lambda) $and the $n$th increment  $\eta^{(n)}(t,\tau)$ admit the spectral
representations
\[\eta(t)=\int_{-\infty}^{\infty}e^{i\lambda t}dZ_{\eta}(\lambda)\] and
\[\eta^{(n)}(t,\tau)=\int_{-\infty}^{\infty}e^{i\lambda t}(1-e^{-i\lambda\tau})^ndZ_{\eta}(\lambda),\]
where $Z_{\eta}(\lambda)$   is a random process with independent increments defined on
$\mr R$ corresponding to the spectral function $G(\lambda)$.
The stationary increment process $\zeta^{(n)}(t,\tau)$ admits the
spectral representation
\[\zeta^{(n)}(t,\tau)=\int_{-\infty}^{\infty}e^{i\lambda t}(1-e^{-i\lambda\tau})^n\frac{(1+i\lambda )^n}{(i\lambda)^n}dZ_{\xi^{(n)}+\eta^{(n)}}(\lambda)\]\[=\int_{-\infty}^{\infty}e^{i\lambda t}(1-e^{-i\lambda\tau})^n\frac{(1+i\lambda )^n}{(i\lambda)^n}dZ_{\xi^{(n)}}(\lambda)+\int_{-\infty}^{\infty}e^{i\lambda t}(1-e^{-i\lambda\tau})^n dZ_{\eta }(\lambda),\]
where $dZ_{\eta^{(n)}}(\lambda)=\dfrac{\lambda^{2n}} {(1+\lambda^2 )^n}dZ_{\eta }(\lambda)$,
$\lambda\in\mr R$. From the described representations we can obtain the spectral density $p(\lambda)$ of the random
process $\zeta(t)$ :
\[p(\lambda)=f(\lambda)+\frac{\lambda^{2n}} {(1+\lambda^2 )^n}g(\lambda).\]
Finally, the  spectral representation of the estimate $\widehat{A}\eta$ is
\[\widehat{A}\eta=\int_{-\infty}^{\infty}
h_{\tau}(\lambda)dZ_{\xi^{(n)}+\eta^{(n)}}(\lambda),\]
and then the  spectral representation of the estimate $\widehat{A}\xi$ of the functional $A\xi$ is
\be \label{otsinka A f_cont_st.n}
\widehat{A}\xi=A\zeta-\int_{-\infty}^{\infty}
h_{\tau}(\lambda)dZ_{\xi^{(n)}+\eta^{(n)}}(\lambda).\ee

It comes from  condition 2) that for all $t\leq0$ the function $h_{\tau}(\lambda)$ satisfies the following equality:
\be\label{int_cont_umova_2}
\int_{-\infty}^{\infty}\ld[A(\lambda)\frac{(-i\lambda)^n}{(1-i\lambda )^n}g(\lambda)-h_{\tau}(\lambda)\ld(f(\lambda)+\frac{\lambda^{2n}} {(1+\lambda^2 )^n}g(\lambda)\rd)\rd]\frac{(1-i\lambda )^n}{(-i\lambda)^{n}}(1-e^{i\lambda\tau})^n e^{-i\lambda
t}d\lambda=0.
\ee
 Let us define a function
 \[C^{\tau}(\lambda)=\ld[A(\lambda)\frac{(-i\lambda)^n}{(1-i\lambda )^n}g(\lambda)-h_{\tau}(\lambda)\ld(f(\lambda)+\frac{\lambda^{2n}} {(1+\lambda^2 )^n}g(\lambda)\rd)\rd]\frac{(1-i\lambda )^n}{(-i\lambda)^{n}}(1-e^{i\lambda\tau})^n, \quad \lambda\in \mr R,\] and its Fourier transformation \[\me c_{\tau}(t)=\int
_{-\infty}^{\infty}C^{\tau}(\lambda)e^{-i\lambda t}d\lambda,\quad
t\in\mr R.\] Condition ($\ref{int_cont_umova_2}$) let us conclude that the function
$\me c_{\tau}(t)$ is equal to $0$ on $(-\infty;0]$, which implies
\[C^{\tau}(\lambda)=\int _{0}^{\infty}\me c_{\tau}(t)e^{i\lambda
t}dt.\]
Hence the function
$h_{\tau}(\lambda)$ has a form
\[h_{\tau}(\lambda)=
 \frac{A(\lambda)(-i\lambda)^ng(\lambda)}{(1-i\lambda )^n\ld(f(\lambda)+\frac{\lambda^{2n}}{(1+\lambda^2 )^n}g(\lambda)\rd)}-\frac{(1-e^{i\lambda\tau})^{-n}(-i\lambda)^{n}C^{\tau}(e^{i\lambda})}{(1-i\lambda )^n
 \ld(f(\lambda)+\frac{\lambda^{2n}}{(1+\lambda^2 )^n}g(\lambda)\rd)},\]
\[A(\lambda)=\int_{0}^{\infty}a(t)e^{-i\lambda t}dt.
\]

Define a closed linear subspace $L_2^{0}\ld(f(\lambda)+\frac{\lambda^{2n}}{{(1+\lambda^2)^n}}g(\lambda)\rd)$
of the Hilbert space
$L_2\ld(f(\lambda)+\frac{\lambda^{2n}}{{(1+\lambda^2)^n}}g(\lambda)\rd)$ generated by the set of functions
\[\ld\{e^{i\lambda t}(1-e^{-i\lambda\tau})^n\dfrac{(1+i\lambda)^n}{(i\lambda)^{n}}:t\leq0\rd\}.\]
Using condition 1) and isometry between the subspaces $L_2^{0}\ld(f(\lambda)+\frac{\lambda^{2n}}{{(1+\lambda^2)^n}}g(\lambda)\rd)$ and $H^{0}(\xi^{(n)}_{\tau}+\eta^{(n)}_{\tau})$  we can obtain
 the following properties of the spectral characteristic
$h_{\tau}(\lambda)$:
\[
h_{\tau}(\lambda)=h(\lambda)(1-e^{-i\lambda\tau})^n\frac{(1+i\lambda )^n}{(i\lambda)^n},\quad
h(\lambda)=\int_{-\infty}^{0}s(t)e^{i\lambda t}dt,
\]
for some function  $s(t)\in L_2^{0}$,
\[
 \int_{-\infty}^{\infty} |h_{\tau}(\lambda)|^2|1-e^{i\lambda\tau}|^{2n}\ld(\frac{(1+\lambda^2 )^n}{\lambda^{2n}}f(\lambda)+g(\lambda)\rd)d\lambda<\infty,
\quad \frac{(i\lambda)^nh_{\tau}(\lambda)}{(1+i\lambda)^{n}(1-e^{-i\lambda
\tau})^n}\in L_2^{0}.
\]
Thus for all $s\geq0$ the following equality holds true
\[ \int_{-\infty}^{\infty} \ld[\frac{A(\lambda)(1-e^{-i\lambda\tau})^{-n}\lambda^{2n}g(\lambda)}{(1+\lambda^2 )^n\ld(f(\lambda)+\frac{\lambda^{2n}}{(1+\lambda^2 )^n}g(\lambda)\rd)}-\frac{|1-e^{i\lambda\tau}|^{-n}\lambda^{2n}C^{\tau}(\lambda)}{(1+\lambda)^{2n}
 \ld(f(\lambda)+\frac{\lambda^{2n}}{(1+\lambda^2)^n}g(\lambda)\rd)}\rd]e^{-i\lambda
s}d\lambda=0.\]
The last condition provides us with an equation defining the  function $\me c_{\tau}(t)$:
\[\int_{-\tau n}^{\infty}
   \me a_{\tau}(t)\int_{-\infty}^{\infty}e^{-i\lambda(t+s)}\frac{\lambda^{2n}g(\lambda)}{|1-e^{i\lambda\tau}|^{2n}
   \ld((1+\lambda^2)^nf(\lambda)+\lambda^{2n}g(\lambda)\rd)}d\lambda
dt\]
\be\label{equation_c(t)_fil_cont_st_n}=\int_{0}^{\infty}\me
c_{\tau}(t)\int_{-\infty}^{\infty}e^{i\lambda(t-s)}\frac{\lambda^{2n}}{|1-e^{i\lambda\tau}|^{2n}
\ld((1+\lambda^2)^n f(\lambda)+\lambda^{2n}g(\lambda)\rd)}d\lambda
dt, \quad s\geq0,\ee
where \be\label{function_a_tau_f_cont_st_n}\me a_{\tau}(t)=\sum_{l=\max\ld\{0,\ld[-\frac{t}{\tau}\rd]'\rd\}}^n C_n^l(-1)^la(t+\tau l),\quad t\geq-\tau n.\ee
Here by $[x]'$ we denote the least integer number among numbers which are
greater than or equal to $x$. Equation $(\ref{equation_c(t)_fil_cont_st_n})$ can be rewritten in terms of linear operators in Hilbert space $L_2$. For functions $\me x(t)\in L_2[-\tau n;\infty)$ and $\me y(t), \me z(t)\in L_2[0;\infty)$ define linear operators
\[(\me S^{\tau}_{\infty}\me x)(s)=\int_{-\tau n}^{\infty}\me
x(t)\int_{-\infty}^{\infty}e^{-i\lambda(t+s)}\frac{\lambda^{2n}g(\lambda)}{|1-e^{i\lambda\tau}|^{2n}
   \ld((1+\lambda^2)^nf(\lambda)+\lambda^{2n}g(\lambda)\rd)}d\lambda
dt,\, s\in[0;\infty),\]
\[(\me P^{\tau}_{\infty}\me y)(s)=\int_{0}^{\infty}\me
y(t)\int_{-\infty}^{\infty}e^{i\lambda(t-s)}\frac{\lambda^{2n}}{|1-e^{i\lambda\tau}|^{2n}
\ld( (1+\lambda^2)^n f(\lambda)+\lambda^{2n}g(\lambda)\rd)}d\lambda
dt,\, s\in[0;\infty),\]
\[(\me Q_{\infty}\me z)(s)=\int_{0}^{\infty}\me
z(t)\int_{-\infty}^{\infty}e^{i\lambda(t-s)}\frac{(1+\lambda^2)^n f(\lambda)g(\lambda)}{
 (1+\lambda^2)^n f(\lambda)+\lambda^{2n}g(\lambda)}d\lambda
dt,\, s\in[0;\infty).\]
Then equation $(\ref{equation_c(t)_fil_cont_st_n})$ can be written as
\[(\me S^{\tau}_{\infty}\me a_{\tau})(s)=(\me P^{\tau}_{\infty}\me c_{\tau})(s), \quad s\geq0,\]
and its solution is
\[\me c_{\tau}(s)=((\me P^{\tau}_{\infty})^{-1}\me S^{\tau}_{\infty}\me a_{\tau})(s), \quad s\geq0.\]

The obtained function $\me c_{\tau}(s)$ allows us to write a solution to filtering problem. The spectral characteristic
$h_{\tau}(\lambda)$ of the optimal estimate $\widehat{A}\xi$ of the functional
$A\xi$ is calculated by the formula
\be \label{int_cont_spectr-char1} h_{\tau}(\lambda)=
 \frac{A(\lambda)(1+i\lambda )^n(-i\lambda)^ng(\lambda)}{(1+\lambda^2 )^n f(\lambda)+ \lambda^{2n} g(\lambda)}-\frac{(1+i\lambda )^n(-i\lambda)^{n}C^{\tau}(e^{i\lambda})}{(1-e^{i\lambda\tau})^n
 \ld((1+\lambda^2 )^n f(\lambda)+ \lambda^{2n} g(\lambda)\rd)}, \ee
where \[C^{\tau}(\lambda)=\int_{0}^{\infty}((\me P^{\tau}_{\infty})^{-1}\me S^{\tau}_{\infty}\me a_{\tau})(t)e^{i\lambda
t}dt.\]
The value of the
mean-square error of the estimate $\widehat{A}\xi$ of the functional $A\xi$ is calculated by the
formula
\[\Delta\ld(f,g;\widehat{A}\xi\rd)=\Delta\ld(f,g;\widehat{A}\eta\rd)= \mt E\ld|A\eta-\widehat{A}\eta\rd|^2
\]
\[=\frac{1}{2\pi}\int_{-\infty}^{\infty}
\frac{\ld|A(\lambda)(1-e^{i\lambda
\tau})^{n}(1+\lambda^2)^nf(\lambda)+\lambda^{2n}C^{\tau}(\lambda)\rd|^2}{|1-e^{i\lambda\tau}|^{2n}(1+\lambda^2)^{2n}
(f(\lambda)+\frac{\lambda^{2n}}{(1+\lambda^2)^n}g(\lambda))^2}g(\lambda)d\lambda
\]\[
+\frac{1}{2\pi}\int_{-\infty}^{\infty}\frac{\ld|A(\lambda)(1-e^{i\lambda
\tau})^{n}(-i\lambda)^{n}g(\lambda)-(-i\lambda)^{n}C^{\tau}(\lambda)\rd|^2}{|1-e^{i\lambda\tau}|^{2n}(1+\lambda^2)^n
(f(\lambda)+\frac{\lambda^{2n}}{(1+\lambda^2)^n}g(\lambda))^2}
f(\lambda)d\lambda\]
\be\label{int_cont_pohybka1}=\langle \me S^{\tau}_{\infty}\me
 a_{\tau},(\me P^{\tau}_{\infty})^{-1}\me S^{\tau}_{\infty}\me a_{\tau}\rangle+\langle\me Q_{\infty}\me a,\me
 a\rangle,\ee
where a function $\me a(t)$, $t\geq0$, is defined as $\me a(t)=a(t)$.

Summing up we can formulate the following theorem.

\begin{thm}\label{int_cont_thm1}
Consider a random process $\{\xi(t),t\in\mr R\}$ with stationary
$n$th increments $\xi^{(n)}(t,\tau)$ which has spectral density
$f(\lambda)$. Let $\{\eta(t),t\in\mr R\}$ be a stationary process with spectral density $g(\lambda)$. Assume that these processes are uncorrelated and their spectral densities satisfy minimality conditions (\ref{int_cont_umova1}). If the function $\me a_{\tau}(t)$ defined by $(\ref{function_a_tau_f_cont_st_n})$ satisfies conditions $(\ref{cond_fro_function_a_f_cont_st_n})$, the optimal linear estimate $\widehat{A}\xi$ of the functional
$A\xi$ of unknown values $\xi(t)$, $t\leq0$,
from observations of a process $\xi(t)+\eta(t)$ at points $t\leq0$ is determined by formula
(\ref{otsinka A f_cont_st.n}). The spectral characteristic
$h_{\tau}(\lambda)$ of the optimal estimate $\widehat{A}\xi$ can be calculated by formula (\ref{int_cont_spectr-char1}). The value of the mean-square error $\Delta(f,g;\widehat{A}\xi)$ of the optimal estimate
is calculated by formula (\ref{int_cont_pohybka1}).
\end{thm}

As a particular case of Theorem $\ref{int_cont_thm1}$ we can derive an estimate
\be \label{otsinka A_T f_cont_st.n}
\widehat{A}_T\xi=A_T\zeta-\int_{-\infty}^{\infty}
h_{\tau,T}(\lambda)dZ_{\xi^{(n)}+\eta^{(n)}}(\lambda)\ee
of the functional $A_T\xi=\int_0^Ta(t)\xi(-t)dt$
from observations of the process
$\xi (t)+\eta (t)$ at  points $t\leq0$. Consider a function $a(t)$, $t\geq0$, which is equal to $0$ for $t>T$. In this case the formula for calculating the spectral characteristic $h_{\tau,T}(\lambda)$ of the
linear estimate is the following:
\be \label{int_cont_spectr-char1_A_T} h_{\tau,T}(\lambda)=
\frac{A_T(\lambda)(1+i\lambda )^n(-i\lambda)^ng(\lambda)}{(1+\lambda^2 )^n f(\lambda)+ \lambda^{2n} g(\lambda)}-\frac{(1+i\lambda )^n(-i\lambda)^{n}C^{\tau}_T(e^{i\lambda})}{(1-e^{i\lambda\tau})^n
 \ld((1+\lambda^2 )^n f(\lambda)+ \lambda^{2n} g(\lambda)\rd)}, \ee
\[A_T(\lambda)=\int_{0}^{T}a(t)e^{-i\lambda t}dt,\quad C^{\tau}_T(\lambda)=\int_{0}^{\infty}((\me P^{\tau}_{\infty})^{-1}\me S^{\tau}_{T}\me a_{\tau,T})(t)e^{i\lambda
t}dt,\]
where a linear operator $S^{\tau}_{T}$ is defined as
\[(\me S^{\tau}_{T}\me x)(s)=\int_{-\tau n}^{T}\me
x(t)\int_{-\infty}^{\infty}e^{-i\lambda(t+s)}\frac{\lambda^{2n}g(\lambda)}{|1-e^{i\lambda\tau}|^{2n}
   \ld((1+\lambda^2)^n f(\lambda)+\lambda^{2n}g(\lambda)\rd)}d\lambda
dt,\, s\in[0;\infty),\]
and a function $\me a_{\tau,T}(t)$, $t\in[-\tau n;T]$, is defined as
\[\me a_{\tau,T}(t)=\sum_{l=\max\ld\{0,\ld[-\frac{t}{\tau}\rd]'\rd\}}^{\min\ld\{n,\ld[\frac{T-t}{\tau}\rd]\rd\}} C_n^l(-1)^la(t+\tau l),\quad t\in[-\tau n;T].\]
The mean-square error of the estimate $\widehat{A}_T\xi$ is calculated by the formula
\[\Delta\ld(f,g;\widehat{A}_T\xi\rd)=\Delta\ld(f,g;\widehat{A}_T\eta\rd)= \mt E\ld|A_T\eta-\widehat{A}_T\eta\rd|^2
\]
\[=\frac{1}{2\pi}\int_{-\infty}^{\infty}
\frac{\ld|A_T(\lambda)(1-e^{i\lambda
\tau})^{n}(1+\lambda^2)^nf(\lambda)+\lambda^{2n}C^{\tau}_T(\lambda)\rd|^2}{|1-e^{i\lambda\tau}|^{2n}(1+\lambda^2)^{2n}
(f(\lambda)+\frac{\lambda^{2n}}{(1+\lambda^2)^n}g(\lambda))^2}g(\lambda)d\lambda
\]\[
+\frac{1}{2\pi}\int_{-\infty}^{\infty}\frac{\ld|A_T(\lambda)(1-e^{i\lambda
\tau})^{n}(-i\lambda)^{n}g(\lambda)-(-i\lambda)^{n}C^{\tau}_T(\lambda)\rd|^2}{|1-e^{i\lambda\tau}|^{2n}(1+\lambda^2)^n
(f(\lambda)+\frac{\lambda^{2n}}{(1+\lambda^2)^n}g(\lambda))^2}
f(\lambda)d\lambda\]
\be\label{int_cont_pohybka1_A_T}=\langle \me S^{\tau}_{T}\me
 a_{\tau},(\me P^{\tau}_{\infty})^{-1}\me S^{\tau}_{T}\me a_{\tau,T}\rangle+\langle\me Q_{T}\me a_T,\me
 a_T\rangle,\ee
where a linear operator $Q_{T}$ is defined as
\[(\me Q_{T}\me z)(s)=\int_{0}^{T}\me
z(t)\int_{-\infty}^{\infty}e^{i\lambda(t-s)}\frac{(1+\lambda^2)^n f(\lambda)g(\lambda)}{(1+\lambda^2)^nf(\lambda)+\lambda^{2n}g(\lambda)}d\lambda
dt,\, s\in[0;\infty),\]
and a function $\me a_{T}(t)$, $t\in[0;T]$, is defined as $\me a_{T}(t)=a(t)$.

The following theorem holds true.

\begin{thm}\label{int_cont_thm2}
Consider a random process $\{\xi(t),t\in\mr R\}$ with stationary
$n$th increments $\xi^{(n)}(t,\tau)$ which has spectral density
$f(\lambda)$. Let $\{\eta(t),t\in\mr R\}$ be a stationary process with spectral density $g(\lambda)$. Assume that these processes are uncorrelated and their spectral densities satisfy minimality conditions (\ref{int_cont_umova1}). Then the optimal linear estimate $\widehat{A}_T\xi$ of the functional
$A_T\xi$ of unknown values $\xi(t)$, $t\in[-T;0]$,
from observations of a process $\xi(t)+\eta(t)$ at points $t\leq0$ is determined by formula
(\ref{otsinka A_T f_cont_st.n}). The spectral characteristic
$h_{\tau,T}(\lambda)$ of the optimal estimate $\widehat{A}_T\xi$ can be calculated by formula (\ref{int_cont_spectr-char1_A_T}). The value of the mean-square error $\Delta(f,g;\widehat{A}_T\xi)$ of the optimal estimate
is calculated by formula (\ref{int_cont_pohybka1_A_T}).
\end{thm}

\section{Minimax-robust method of filtering}

In the previous section we solved filtering problem for the functionals
${A}\xi$ and ${A}_T\xi$ of unknown values of the random
process  $\xi(t)$ based on observations of the random process
$\xi(t)+\eta(t)$. Derived formulas $(\ref{int_cont_spectr-char1_A})$
and $(\ref{int_cont_spectr-char1_A_T})$ for spectral characteristics of the mean-square estimates
$\widehat{A}\xi$ and $\widehat{A}_T\xi$ can be used only if we know spectral functions $f(\lambda)$ and $g(\lambda)$ of the random processes $\xi(t)$ and $\eta(t)$. In practice we do not the exact view of these spectral densities. Therefore in this section we propose the minimax (robust) approach to
estimation of functionals which depend on the unknown values of
random process with stationary increments. This method allows us to find an estimate that minimizes
the maximum of mean-square errors of the estimates for all spectral densities from a
given class $\md D=\md D_f\times\md D_g$  of admissible spectral
densities.

\begin{ozn} For a given class of spectral densities $\mathcal{D}=\md
D_f\times\md D_g$  spectral densities
$f^0(\lambda)\in\mathcal{D}_f$ and $g^0(\lambda)\in\md D_g$ are called
least favorable in the class $\mathcal{D}$ for the optimal linear
filtering of the functional $A\xi$ if the following relation holds
true
\[\Delta(f^0,g^0)=\Delta\ld(h(f^0,g^0);f^0,g^0\rd)=
\max_{(f,g)\in\mathcal{D}_f\times\md
D_g}\Delta\ld(h(f,g);f,g\rd).\]\end{ozn}

\begin{ozn}
For a given class of spectral densities $\mathcal{D}=\md
D_f\times\md D_g$ the spectral characteristic $h^0(\lambda)$ of
the optimal linear estimate of the functional $A\xi$ is called
minimax-robust if there are satisfied conditions
\[h^0(\lambda)\in H_{\mathcal{D}}=\bigcap_{(f,g)\in\mathcal{D}_f\times\md D_g}L_2^{0}\ld(f(\lambda)+\frac{\lambda^{2n}}{(1+\lambda^2)^n}g(\lambda)\rd),\]
\[\min_{h\in H_{\mathcal{D}}}\max_{(f,g)\in \mathcal{D}_f\times\md D_g}\Delta(h;f,g)=\max_{(f,g)\in\mathcal{D}_f\times\md
D_g}\Delta(h^0;f,g).\]\end{ozn}

From the formulas proposed in the previous section and the
introduced definitions we can obtain the following statement.

\begin{lem}Spectral densities $f^0\in\mathcal{D}_f$,
$g^0\in\mathcal{D}_g$ satisfying minimality condition $(\ref{int_cont_umova1})$
are least favorable in the class $\md D=\md D_f\times\md D_g$ for
the optimal linear filtering of the functional $A\xi$ if operators
$ (\me P^{\tau}_{\infty})^0$, $(\me S^{\tau}_{\infty})^0$, $(\me Q_{\infty})^0$  determined by the
Fourier coefficients of the functions
\[\dfrac{\lambda^{2n}|1-e^{i\lambda\tau}|^{-2n}}{
 (1+\lambda^2)^n f^0(\lambda)+\lambda^{2n}g^0(\lambda)}, \quad
\dfrac{\lambda^{2n}g^0(\lambda)|1-e^{i\lambda\tau}|^{-2n}}
{(1+\lambda^2)^nf^0(\lambda)+\lambda^{2n}g^0(\lambda)},\quad
\dfrac{(1+\lambda^2)^n f^0(\lambda)g^0(\lambda)}{
 (1+\lambda^2)^n f^0(\lambda)+\lambda^{2n}g^0(\lambda)}\]
determine a solution of the conditional extremum problem
\be\max_{(f,g)\in \mathcal{D}_f\times\md D_g}\ld(\langle \me S^{\tau}_{\infty}\me
 a_{\tau},(\me P^{\tau}_{\infty})^{-1}\me S^{\tau}_{\infty}\me a_{\tau}\rangle+\langle\me Q_{\infty}\me a,\me
 a\rangle\rd)= \langle (\me S^{\tau}_{\infty})^0\me
 a_{\tau},((\me P^{\tau}_{\infty})^0)^{-1}(\me S^{\tau}_{\infty})^0\me a_{\tau}\rangle+\langle(\me Q_{\infty})^0\me a,\me
 a\rangle. \label{minimax1 fs}\ee The minimax spectral
characteristic is determined as $h^0=h_{\tau}(f^0,g^0)$ if
$h_{\tau}(f^0,g^0)\in H_{\mathcal{D}}$.
\end{lem}

The minimax-robust spectral characteristic $h^0$ and the least favorable spectral densities $(f^0,g^0)$ form a saddle point of
the function $\Delta(h;f,g)$ on the set
$H_{\mathcal{D}}\times\mathcal{D}$.
The saddle point inequalities
\[\Delta(h;f^0,g^0)\geq\Delta(h^0;f^0,g^0)\geq\Delta(h^0;f,g)
\quad\forall f\in \mathcal{D}_f,\forall g\in \mathcal{D}_g,\forall
h\in H_{\mathcal{D}}\] hold true if $h^0=h_{\tau}(f^0,g^0)$ and
$h_{\tau}(f^0,g^0)\in H_{\mathcal{D}}$, where a pair $(f^0,g^0)$  provides a
solution to the following conditional extremum problem:

\[\widetilde{\Delta}(f,g)=-\Delta\ld(h_{\tau}(f^0,g^0);f,g\rd)\to
\inf,\quad (f,g)\in \mathcal{D}_f\times \mathcal{D}_g,\]
\[
\Delta\ld(h_{\tau}(f^0,g^0);f,g\rd)\]\[=\frac{1}{2\pi}\int_{-\infty}^{\infty}
\frac{\ld|A(\lambda)(1-e^{i\lambda
\tau})^{n}(1+\lambda^2)^nf^0(\lambda)+\lambda^{2n}\int_{0}^{\infty}(((\me P^{\tau}_{\infty})^0)^{-1}(\me S^{\tau}_{\infty})^0\me a_{\tau})(t)e^{i\lambda
t}dt\rd|^2}{|1-e^{i\lambda\tau}|^{2n}(1+\lambda^2)^{2n}
(f^0(\lambda)+\frac{\lambda^{2n}}{(1+\lambda^2)^n}g^0(\lambda))^2}g(\lambda)d\lambda
\]\[
+\frac{1}{2\pi}\int_{-\infty}^{\infty}\frac{\ld|A(\lambda)(1-e^{i\lambda
\tau})^{n}(-i\lambda)^{n}g^0(\lambda)-(-i\lambda)^{n}\int_{0}^{\infty}(((\me P^{\tau}_{\infty})^0)^{-1}(\me S^{\tau}_{\infty})^0\me a_{\tau})(t)e^{i\lambda
t}dt\rd|^2}{|1-e^{i\lambda\tau}|^{2n}(1+\lambda^2)^n
(f^0(\lambda)+\frac{\lambda^{2n}}{(1+\lambda^2)^n}g^0(\lambda))^2}
f(\lambda)d\lambda.\]

This last extremum problem is equivalent to the unconditional
extremum problem
\[\Delta_{\mathcal{D}}(f,g)=\widetilde{\Delta}(f,g)+ \delta(f,g|\mathcal{D}_f\times
\mathcal{D}_g)\to\inf,\] where $\delta(f,g|\mathcal{D}_f\times
\mathcal{D}_g)$ is the indicator function of the set
$\mathcal{D}_f\times\mathcal{D}_g$.
 Solution $(f^0,g^0)$ to this unconditional
extremum problem is characterized by condition   $0\in
\partial\Delta_{\mathcal{D}}(f^0,g^0)$. For more details we refer to the books \cite{Moklyachuk:2008nonsm,Pshenychn}.

\

\section{Least favorable spectral densities in the class $\md D^0_f\times\md D^0_g$}

Consider a  set of admissible spectral densities $\md D=\md
D^0_f\times\md D^0_g$, where
\[\md D^0_f=\ld\{f(\lambda)|\frac{1}{2\pi}\int_{-\infty}^{\infty}
f(\lambda)d\lambda\leq P_1\rd\},\quad\md
D^0_g=\ld\{g(\lambda)|\frac{1}{2\pi}\int_{-\infty}^{\infty} g(\lambda)d\lambda\leq
P_2\rd\}.\] Assume that spectral densities $f^0\in\md D^0_f$,
$g^0\in\md D^0_g$ and functions \be
h_{\tau,f}(f^0,g^0)=\frac{\ld|A(\lambda)(1-e^{i\lambda
\tau})^{n}(-i\lambda)^{n}g^0(\lambda)-(-i\lambda)^{n}\int_{0}^{\infty}(((\me P^{\tau}_{\infty})^0)^{-1}(\me S^{\tau}_{\infty})^0\me a_{\tau})(t)e^{i\lambda
t}dt\rd|}{|1-e^{i\lambda\tau}|^{n}|1+i\lambda|^n
(f^0(\lambda)+\frac{\lambda^{2n}}{(1+\lambda^2)^n}g^0(\lambda))},\label{hf_interp_shum
fs}\ee
 \be h_{\tau,g}(f^0,g^0)=\frac{\ld|A(\lambda)(1-e^{i\lambda
\tau})^{n}(1+\lambda^2)^nf^0(\lambda)+\lambda^{2n}\int_{0}^{\infty}(((\me P^{\tau}_{\infty})^0)^{-1}(\me S^{\tau}_{\infty})^0\me a_{\tau})(t)e^{i\lambda
t}dt\rd|}{|1-e^{i\lambda\tau}|^{n}(1+\lambda^2)^{n}
(f^0(\lambda)+\frac{\lambda^{2n}}{(1+\lambda^2)^n}g^0(\lambda))}\label{hg_interp_shum
fs}\ee are bounded. Under these conditions the functional
$\Delta(h_{\tau}(f^0,g^0);f,g)$ is continuous and bounded in $\md
L_1\times\md L_1$ space. The condition
$0\in\partial\Delta_{\md D}(f^0,g^0)$ implies that least favorable
densities $f^0\in\md D^0_f$, $g^0\in\md D^0_g$ satisfy the equation
\[\ld|A(\lambda)(1-e^{i\lambda
\tau})^{n}(1+\lambda^2)^nf^0(\lambda)+\lambda^{2n}\int_{0}^{\infty}(((\me P^{\tau}_{\infty})^0)^{-1}(\me S^{\tau}_{\infty})^0\me a_{\tau})(t)e^{i\lambda
t}dt\rd|\]
\be =\alpha_1|1-e^{i\lambda\tau}|^{n}
\ld((1+\lambda^2)^n f^0(\lambda)+\lambda^{2n}g^0(\lambda)\rd),\label{D1
rivn1 fs}\ee
\[\ld|A(\lambda)(1-e^{i\lambda
\tau})^{n}(-i\lambda)^{n}g^0(\lambda)-(-i\lambda)^{n}\int_{0}^{\infty}(((\me P^{\tau}_{\infty})^0)^{-1}(\me S^{\tau}_{\infty})^0\me a_{\tau})(t)e^{i\lambda
t}dt\rd|\] \be
=\alpha_2|1-e^{i\lambda\tau}|^{n}|1-i\lambda|^{-n}
\ld((1+\lambda^2)^n f^0(\lambda)+\lambda^{2n}g^0(\lambda)\rd),\label{D1
rivn2 fs}\ee where $\alpha_1\geq0$ and $\alpha_2\geq0$ are constants
such that
 $\alpha_1\neq0$ if $(2\pi)^{-1}\int_{-\infty}^{\infty} f^0(\lambda)d\lambda=P_1$ and $\alpha_2\neq0$
if $(2\pi)^{-1}\int_{-\infty}^{\infty} g^0(\lambda)d\lambda=P_2$.

Thus, the following statements hold true.

\begin{thm} Let conditions ($\ref{cond_fro_function_a_f_cont_st_n}$) hold true, spectral densities $f^0(\lambda)\in\md D^0_f$,
$g^0(\lambda)\in\md D^0_g$ satisfy  conditions $(\ref{int_cont_umova1})$ and functions $h_{\tau,f}(f^0,g^0)$,
$h_{\tau,g}(f^0,g^0)$ be bounded. The spectral densities
$f^0(\lambda)$ and $g^0(\lambda)$ are  least favorable in
the class $\md D=\md D^0_f \times\md D^0_g$ for the optimal linear
estimation of the functional $A\xi$ if they are a solution of equations
$(\ref{D1 rivn1 fs})$-$(\ref{D1 rivn2 fs})$ and determine a solution to
extremum problem $(\ref{minimax1 fs})$. Minimax spectral characteristic $h_{\tau}(f^0,g^0)$ of the optimal estimate of the
functional $A\xi$ is determined by formula $(\ref{int_cont_spectr-char1})$.\end{thm}

\begin{thm}
Let known spectral density $f(\lambda) $ and spectral density
$g^0(\lambda)\in\md D^0_g$ satisfy conditions $(\ref{int_cont_umova1})$. Let conditions ($\ref{cond_fro_function_a_f_cont_st_n}$) hold true and the function $h_{\tau,g}(f ,g^0)$ be bounded.
Then spectral density $g^0(\lambda)$ is least favorable in the class
 $\md D^0_g$ for the optimal linear filtering of the  functional $A\xi$
 if it has the form
\[g^0(\lambda)=(1+\lambda^2)^n\lambda^{-2n}\max\ld\{0,f_1(\lambda)\rd\},\]
\be\label{minimax_dens_f1}f_1(\lambda)=\frac{\ld|A(\lambda)(1-e^{i\lambda
\tau})^{n}(1+\lambda^2)^nf(\lambda)+\lambda^{2n}\int_{0}^{\infty}(((\me P^{\tau}_{\infty})^0)^{-1}(\me S^{\tau}_{\infty})^0\me a_{\tau})(t)e^{i\lambda
t}dt\rd|}{\alpha_2|1-e^{i\lambda\tau}|^{n}(1+\lambda^2)^n}-f(\lambda)\ee
and the
pair $(f,g^0)$ determines a solution to extremum problem
$(\ref{minimax1 fs})$. Minimax spectral characteristic $h_{\tau}(f,g^0)$ of the optimal estimate of the
functional $A\xi$ is determined by formula $(\ref{int_cont_spectr-char1})$.\end{thm}

\begin{thm}Let spectral density
$f^0(\lambda)\in\md D^0_f$ and known spectral density $g(\lambda) $ satisfy conditions $(\ref{int_cont_umova1})$. Let conditions ($\ref{cond_fro_function_a_f_cont_st_n}$) hold true and the function $h_{\tau,f}(f^0,g)$ be bounded. Then
spectral density $f^0(\lambda)$ is least favorable in the class
 $\md D^0_f$ for the optimal linear filtering of the  functional $A\xi$
 if it has the form
\[f^0(\lambda)=(1+\lambda^2)^n\lambda^{-2n}\max\ld\{0,g_1(\lambda)\rd\},\]
\be\label{minimax_dens_g1}g_1(\lambda)=\frac{\ld|A(\lambda)(1-e^{i\lambda
\tau})^{n}(-i\lambda)^{n}g(\lambda)-(-i\lambda)^{n}\int_{0}^{\infty}(((\me P^{\tau}_{\infty})^0)^{-1}(\me S^{\tau}_{\infty})^0\me a_{\tau})(t)e^{i\lambda
t}dt\rd|}{\alpha_1|1-e^{i\lambda\tau}|^{n}|1+i\lambda|^n}-g
(\lambda)\ee
 and the pair $(f^0 ,g)$ determines a solution to extremum problem $(\ref{minimax1 fs})$. Minimax spectral characteristic $h_{\tau}(f^0,g)$ of the optimal estimate of the
functional $A\xi$ is determined by formula $(\ref{int_cont_spectr-char1})$.\end{thm}

\section{Least favorable densities in the class $\md D=\md D_{u}^v\times\md D_{\eps}$}

Consider an optimal linear filtering of the functional
$A\xi$ in case of the set $\md D=\md D_u^v\times\md
D_{\eps}$ of admissible spectral densities, where
 \[\md
D_u^v=\ld\{f(\lambda)|v(\lambda)\leq f(\lambda)\leq
u(\lambda),\,\frac{1}{2\pi}\int_{-\infty}^{\infty} f(\lambda)d\lambda\leq P_1\rd\},\]
\[\md
D_{\eps}=\ld\{g(\lambda)|g(\lambda)=(1-\eps)g_2(\lambda)+\eps
w(\lambda),\frac{1}{2\pi}\int_{-\infty}^{\infty} g(\lambda)d\lambda\leq P_2\rd\}.\]
Spectral densities $u(\lambda)$, $v(\lambda)$ and $g_2(\lambda)$ from the definition os the set $\md D=\md D_u^v\times\md
D_{\eps}$
are known, fixed and in addition spectral densities $u(\lambda)$ and $v(\lambda)$ are bounded.

Let functions $h_{\tau,f}(f^0,g^0)$ and $h_{\tau,g}(f^0,g^0)$
determined by formulas $(\ref{hf_interp_shum fs})$,
$(\ref{hg_interp_shum fs})$ are bounded for spectral densities $f^0\in \md D_{u}^v$, $g^0\in \md D_{\eps}$. Condition
$0\in\partial\Delta_{\md D}(f^0,g^0)$ imply
equations that define least favorable spectral densities
\[\ld|A(\lambda)(1-e^{i\lambda
\tau})^{n}(1+\lambda^2)^nf^0(\lambda)+\lambda^{2n}\int_{0}^{\infty}(((\me P^{\tau}_{\infty})^0)^{-1}(\me S^{\tau}_{\infty})^0\me a_{\tau})(t)e^{i\lambda
t}dt\rd|\]
\be = |1-e^{i\lambda\tau}|^{n}
\ld((1+\lambda^2)^n f^0(\lambda)+\lambda^{2n}g^0(\lambda)\rd)(\gamma_1(\lambda)+\gamma_2(\lambda)+\alpha_1),\label{D2
rivn1 fs}\ee
\[\ld|A(\lambda)(1-e^{i\lambda
\tau})^{n}(-i\lambda)^{n}g^0(\lambda)-(-i\lambda)^{n}\int_{0}^{\infty}(((\me P^{\tau}_{\infty})^0)^{-1}(\me S^{\tau}_{\infty})^0\me a_{\tau})(t)e^{i\lambda
t}dt\rd|\] \be
= |1-e^{i\lambda\tau}|^{n}|1-i\lambda|^{-n}
\ld((1+\lambda^2)^n f^0(\lambda)+\lambda^{2n}g^0(\lambda)\rd)(\varphi(\lambda)+\alpha_2),\label{D2
rivn2 fs}\ee
where
$\gamma_1\leq0$ and $\gamma_1=0$ if $f^0(\lambda)\geq
v(\lambda)$; $\gamma_2(\lambda)\geq0$ and $\gamma_2=0$ if
$f^0(\lambda)\leq u(\lambda)$; $\varphi(\lambda)\leq0$ and
$\varphi(\lambda)=0$ when $g^0(\lambda)\geq (1-\eps)g_1(\lambda)$.

The following statements hold true.

\begin{thm} Let conditions ($\ref{cond_fro_function_a_f_cont_st_n}$) hold true, spectral densities $f^0(\lambda)\in \md D_{u}^v$, $g^0(\lambda)\in \md
D_{\eps}$ satisfy  condition $(\ref{int_cont_umova1})$ and functions $h_{\tau,f}(f^0,g^0)$,
$h_{\tau,g}(f^0,g^0)$ determined by formulas
$(\ref{hf_interp_shum fs})$, $(\ref{hg_interp_shum fs})$ be bounded.
Spectral densities $f^0(\lambda)$ and $g^0(\lambda)$ are least
favorable in the class $\md D=\md D_u^v \times\md D_{\eps}$ for the
optimal linear filtering of the functional $A\xi$ if they satisfy
equations $(\ref{D2 rivn1 fs})$, $(\ref{D2 rivn2 fs})$ and determine
a solution to extremum problem $(\ref{minimax1 fs})$. Minimax spectral characteristic $h_{\tau}(f^0,g^0)$ of the optimal estimate of the
functional $A\xi$ is determined by formula $(\ref{int_cont_spectr-char1})$. \end{thm}

\begin{thm} Let known spectral density $f(\lambda) $ and spectral density
$g^0(\lambda)\in\md D_{\eps}$ satisfy conditions $(\ref{int_cont_umova1})$. Let conditions ($\ref{cond_fro_function_a_f_cont_st_n}$) hold true and the function $h_{\tau,g}(f ,g^0)$ be bounded.
Spectral density
$g^0(\lambda)$ is  least favorable in the class
  $\md D_{\eps}$ for the optimal linear filtering of the functional $A\xi$
  if it has the form
\[g^0(\lambda)=\max\ld\{(1-\eps)g_2(\lambda),(1+\lambda^2)^n\lambda^{-2n}f_1(\lambda)\rd\},\]
where function $f_1(\lambda)$ is defined by formula ($\ref{minimax_dens_f1}$), and the pair
$(f ,g^0)$ determines a solution to extremum problem $(\ref{minimax1
fs})$. Minimax spectral characteristic $h_{\tau}(f,g^0)$ of the optimal estimate of the
functional $A\xi$ is determined by formula $(\ref{int_cont_spectr-char1})$.\end{thm}

\begin{thm}
Let spectral density
$f^0(\lambda)\in\md D_u^v$ and known spectral density $g(\lambda) $ satisfy conditions $(\ref{int_cont_umova1})$. Let conditions ($\ref{cond_fro_function_a_f_cont_st_n}$) hold true and the function $h_{\tau,f}(f^0,g)$ be bounded.
Spectral density $f^0(\lambda)$ is  least favorable in the class
 $\md D_u^v$ for the optimal linear filtering of the  functional $A\xi$  if it has the form
 \[f^0(\lambda)=\min\ld\{v(\lambda),\max\ld\{u(\lambda),(1+\lambda^2)^n\lambda^{-2n}g_1(\lambda)\rd\}\rd\},\]
where function $g_1(\lambda)$ is defined by formula ($\ref{minimax_dens_g1}$),
and the pair $(f^0,g)$ determines a solution to  extremum problem
$(\ref{minimax1 fs})$. Minimax spectral characteristic $h_{\tau}(f^0,g)$ of the optimal estimate of the
functional $A\xi$ is determined by formula $(\ref{int_cont_spectr-char1})$.\end{thm}

\section{Conclusions}

In this paper we derived two methods of solving a problem of optimal mean-square estimation of  linear functionals $A{\xi}=\int_{0}^{\infty}a(t)\xi(-t)dt$ and $A_T{\xi}=\int_{0}^Ta(t)\xi(-t)dt$
of  unknown values of a random process $\xi(t)$ with $n$th stationary increments from the observation of a process $\xi(t)+\eta(t)$ at points
$t\leq0$. The first method provides us with formulas for calculating the values of mean-square errors and the spectral characteristics of the estimates of the functionals $A{\xi}$ and $A_T{\xi}$ when  spectral densities of the processes are given. The second one, called minimax-robust method, let us solve the problem in the case when spectral densities are not known, but a set of admissible spectral densities is available.

\end{document}